\documentclass{article}

\usepackage{arxiv}

\usepackage[utf8]{inputenc} 
\usepackage[T1]{fontenc}    
\usepackage{hyperref}       
\usepackage{url}            
\usepackage{booktabs}       
\usepackage{nicefrac}       
\usepackage{microtype}      
\usepackage{lipsum}
\usepackage{graphicx}
\usepackage[ruled]{algorithm2e}
\RequirePackage[OT1]{fontenc}
\RequirePackage{amsthm,amsmath}
\RequirePackage[numbers]{natbib}
\usepackage{wrapfig}
\usepackage{subcaption}

\usepackage{amsmath, amsthm, amsfonts}
\usepackage{amssymb,mathrsfs}
\usepackage{float, multirow, url}
\usepackage{array, hhline}
\usepackage{caption}
\usepackage{sidecap}
\usepackage{chngcntr}
\counterwithout{equation}{section}

\def\be{\begin{equation}}
\def\ee{\end{equation}}
\newtheorem{theorem}{Theorem}
\newtheorem{lemma}{Lemma}

\theoremstyle{definition}

\newtheorem*{lemma*}{Lemma}

\def\v{{\|  }}
\def\l{{\vert }}

\def\bOmg{\boldsymbol{\Omega}}
\def\bgam{\boldsymbol{\Gamma}}
\def\bX{\mathbf{X}}
\def\bY{\mathbf{Y}}
\def\bW{\mathbf{W}}
\def\bZ{\mathbf{Z}}
\def\bS{\mathbf{S}}
\def\bA{\mathbf{A}}
\def\bB{\mathbf{B}}

\def\bI{\mathbf{I}}

\def\bm{\mathbf{m}}
\def\bN{\mathbf{N}}

\def\ba{\mathbf{a}}
\def\bb{\mathbf{b}}
\def\bx{\mathbf{x}}

\def\bE{\mathbf{e}}
\def\bU{\mathbf{U}}

\def\bg{\mathbf{g}}
\def\bw{\mathbf{w}}

\def\bmu{\boldsymbol{\mu}}

\def\bsig{\boldsymbol{\Sigma}}

\def\bbeta{\boldsymbol{\beta}}

\def\bsig{\boldsymbol{\Sigma}}

\def\bSXn{\boldsymbol{S}_{\boldsymbol{X},n}}
\def\bSYn{\boldsymbol{S}_{\boldsymbol{Y},n}}
\def\bSXni{\boldsymbol{S}_{\boldsymbol{X},ni}}
\def\bSXnj{\boldsymbol{S}_{\boldsymbol{X},nj}}
\def\bWXni{\boldsymbol{W}_{\boldsymbol{X},ni}}

\def\bVXn{\boldsymbol{V}_{\boldsymbol{X},n}}
\def\bVYn{\boldsymbol{V}_{\boldsymbol{Y},n}}
\def\bVXni{\boldsymbol{V}_{\boldsymbol{X},ni}}

\def\bWYn{\boldsymbol{W}_{\boldsymbol{Y},n}}
\def\bS{\boldsymbol{S}}

\def\xbi{\bX_{\bB_i}}
\def\x1bi{\bX_{1\bB_i}}
\def\xkbi{\bX_{k\bB_i}}
\def\xnbi{\bX_{n\bB_i}}

\def\bgam{\boldsymbol{\Gamma}}
\def\bmu{\boldsymbol{\mu}}
\newcommand{\E}{\mathrm{E}}

\newcommand{\Var}{\mathrm{Var}}

\newcommand{\tr}{\mathrm{tr}}

\DeclareMathOperator*{\argmin}{arg\,min}

\graphicspath{ {./images/} }

\title{High-dimensional Precision Matrix Estimation with a Known
Graphical Structure}

\author{
 
  Thien-Minh Le\\
  Department of Biostatistics\\
  Harvard T.H. Chan School of Public Health\\
  Boston, Massachusetts, U.S.A.\\
  \texttt{thle@hsph.harvard.edu} \\
   \And
 Ping-Shou Zhong\\
 Department of Mathematics, Statistics and Computer Science\\
  University of Illinois at Chicago\\
  Chicago, Illinois, U.S.A.\\
  \texttt{pszhong@uic.edu} \\
  
}

\begin{document}
\maketitle
\begin{abstract}
A precision matrix is the inverse of a covariance matrix.
In this paper, we study the problem of estimating the precision matrix with a known graphical structure under high-dimensional settings. We propose a simple estimator of the precision matrix based on the connection between the known graphical structure and the precision matrix. We obtain the rates of convergence of the proposed estimators and derive the asymptotic normality of the proposed estimator in the high-dimensional setting when the data dimension grows with the sample size. Numerical simulations are conducted to demonstrate the performance of the proposed method. We also show that the proposed method outperforms some existing methods that do not utilize the graphical structure information.

\end{abstract}


\section{Introduction}\label{sec1}

The precision matrix is the inverse of the covariance matrix, widely used in many statistical procedures such as linear discriminant analysis, hypothesis testing for multivariate mean vectors, and confidence regions for mean vectors. In classical asymptotic settings, the number of observations $n$ is assumed to grow to infinity and much bigger than the number of variables $p$, which is 
assumed to be fixed. The inverse of sample covariance is a biased but consistent estimator of the precision matrix (Bai and Shi, 2011). However, in high-dimensional 
settings, where the number of observations $n$ is often smaller than the number of variables $p$, the sample covariance matrix is not a consistent estimator of the population covariance (Bai and Silverstein, 2010; Johnstone, 2001), 
and the precision matrix can not be estimated using the inverse of the sample covariance matrix due to the singularity of the sample covariance matrix.  
To avoid the singularity issue, the generalized inverse of the sample covariance may be used to estimate the precision matrix. However, its application in high-dimensional settings is very limited because the estimator is not unique, singular, and not sparse. These challenges make the problem of estimation of a precision matrix one of the most fundamental issues in high-dimensional 
statistics. Various estimators have been proposed in the literature. One often-used technique is based on regularization. The shrinkage method is one type of regularization methods. Ledoit and Wolf (2003; 
2004) proposed a shrinkage estimator of covariance matrix via minimizing Frobenius norm of the distance between an underlying population covariance and a weighted average of sample covariance and identity matrix. A direct and optimal shrinkage estimator of the precision matrix is proposed by Bodnar, Gupta, and Parolya (2016). No additional assumption is required on covariance or precision matrix in shrinkage-type estimators.  Another popular regularization approach assumes some sparse structure in the precision matrix. Some regularization approaches estimate the entire precision matrix simultaneously by employing the $L_1$-penalized maximum likelihood (see Banerjees et al. (2008), Yuan and Lin (2007), Friedman et al. (2008) and Lam and Fan (2009)). 
Other methods apply the $L_1$-penalized regularization methods in a column by column fashion, such as CLIME (Constrained L1-minimization for Inverse Matrix Estimation) by Cai et al. (2011), 
SCIO (Sparse Column-wise Inverse Operator) by Liu and Lou (2015), or TIGER (Tuning-Insensitive Graph Estimation and Regression) by Liu and Wang (2017). 
However, studying statistical inference problems based on these regularized estimators is challenging due to the non-smoothness of the $L_1$-penalized regularization methods and the use of tuning parameters.

The methods mentioned above are developed for a general precision matrix and do not take advantage of prior knowledge or information available. For instance, in many genetics and genomics research, some well-known pathway databases, such as KEGG, Reactome, and BioCarta, provide prior knowledge about the 
connections among genes or SNPs, which forms graphical structures that could be considered as prior information for any type of statistical analysis. 
Recently, there is a growing interest in using the known graphical structure to improve the accuracy of statistical estimation and inference.
For example, Li and Li (2008) proposed a variable selection procedure incorporating the known network structure. Some other methods incorporating prior knowledge of network structure in the statistical procedure can be found in Li and Li (2010), Zhou and Song (2016), or Gao et al. (2019). The results obtained in these papers are encouraging and demonstrate that utilizing prior knowledge about 
the network structure does help improve the accuracy in estimation and prediction.

To incorporate the prior knowledge about network structure, we develop an estimator of the precision matrix with a known graphical structure. This type of problem did not gain much attention in the past. 
Under the known network structure assumption, Hastie et al. (2008) (Section 17.3.1) proposed an iterative algorithm to estimate the precision matrix and its inverse by maximizing
the Gaussian likelihood function for the precision matrix with constraints. 
Zhou et al. (2011) also proposed another approach to solve the problem where the precision matrix is estimated by minimizing the likelihood function with a constraint on the zero patterns of the precision matrix 
so that it is consistent with the known graphical structure. Both approaches mentioned above are free of tuning parameters. However, they are based on the likelihood function of multivariate Gaussian distribution 
and need some iterative algorithms to obtain the solutions, and an explicit form of the estimator is not available. Due to the nature of these estimators, asymptotic distributions and inference are not easy to establish.

This paper proposes a new and simple method to estimate the precision matrix when its graphical structure is known. Our approach is straightforward, and an explicit form of the estimator is given. The estimator is constructed without using any likelihood function information. The proposed method can be adapted well in a non-parametric framework whenever a connection between the graphical structure and its precision matrix is available. With our simple estimator, we establish the rate of convergence and asymptotic normality of the proposed estimator. Simulation studies demonstrate that the proposed method outperforms some existing methods without utilizing the prior knowledge of network structure. 

The paper is organized as follows. In Section 2, we introduce the basic notation and propose the new estimating procedure. 
Sections 3 studies the rate of convergence and asymptotic normality of the proposed estimator. Simulation studies are investigated in Section 4. 
Section 5 concludes our work. Complete details of the proofs are in the Appendix Section.

\setcounter{equation}{0} 

\section{Estimators of precision matrix with a known structure}
We first define some notations. 
For any vector $\ba = (a_1, a_2,...,a_p)'\in \mathbb{R}^p$, we denote the $L_k$ norm of the vector $\ba$ as $\v \ba \v_k = (\sum_{i=1}^p \l a_i\l ^k)^{1/k}$, where $p$ and $k$ are  positive integers, 
$ p, k \in \mathbb{Z}^+$.
For convenience, we denote $\v \ba \v_2 $ as $\v \ba \v $. The number of non-zero elements of a $p$-dimensional vector $\ba=(a_1,\cdots, a_p)'$ is denoted as $\v \ba \v_0 = \sum_{i=1}^p \l a_i\l ^0$ where $0^0 = 0$.
For any matrix $\bA_{m \times n} $, $\bA^T $ denotes the transpose of $\bA$, $\v \bA \v_1= \max_{j=1,\cdots,n}\sum_{i=1}^m \l a_{ij}\l$, $\v \bA \v_\infty= \v \bA^T \v_1$, $\l \bA \l_\infty= \max_{i=1,\cdots,m;j=1,\cdots,n} \l a_{ij}\l$,
and $\v \bA \v =\sup_{\v \bx \v =1} \v \bA \bx \v$ where $\bx \in \mathbb{R}^{n \times 1}$.  For two sequences of random variables $\{ X_n \}$ and $\{ Y_n \}$, we denote $X_n = o_p(Y_n)$ if  ${X_n}/{Y_n} = o_p(1)$ and $X_n = O_p(Y_n)$  
if ${X_n}/{Y_n} = O_p(1)$. For sequences of numbers $\{ a_n\}, \{ b_n \}$, we denote $ a_n \asymp b_n$ if there exist positive constants $c_1, c_2 > 0,$ such that $c_1 < \l{a_n}/{b_n} \l < c_2$. For notation simplicity, a generic constant $C$  is used, whose value may vary from line to line throughout the paper.

Define $\bX_k = (X_{k,1},X_{k,2},\cdots,X_{k,p})'$ for $k = 1,2,\cdots, n$ and $\bX=(X_1,X_2,\cdots,X_p)'$. Let $\bX_1,\cdots,\bX_n$ be $n$ observed independent and identically distributed (IID) 
realizations of $\bX$ from a multivariate normal distribution with mean vector $\mu$ and covariance matrix $\bsig\in \mathbb{R}^{p \times p}$, $\bsig > 0$ and $\bOmg = \bsig^{-1}$ is the corresponding precision matrix. 
The goal is to estimate $\bOmg$ with a network structure information $\bA=\bA_0$ being given. Here the network information is given by an adjacent matrix which contains zeros and ones. The 
adjacent matrix $\bA=(a_{ij})$ is a $p\times p$ matrix which describes the connectivities among nodes such that $a_{ij}=1$ if nodes $i$ and $j$ are connected; otherwise $a_{ij}=0$.  For instance, the edge between nodes $i$ and $j$ of the Gaussian graphical model is absent if and only if the $(i,j)$-th position of precision matrix $\bOmg$ is 0 (Lauritzen, 1996). Thus, the specification of 
network structure is equivalent to the specification of the structure of a precision matrix. We will make use of the structure of the precision matrix to estimate $\bOmg$.
Assume that $\bOmg$ is a sparse precision matrix with nonzero positions specified, and for $j = 1,2,...,p$, the total of non zero elements in the $j$-th column is $s_j$ where $s_j$ is at the same order as $s_0$, $s_0 =o(\sqrt{n})$. 
Here $s_0$ represents the order of the numbers of non-zeros in every column of precision matrix.  Let $\bw_i$ denote the $i$-th column of $\bOmg$ so that $\bOmg=(\bw_1,\cdots, \bw_p)$. 
We will estimate $\bOmg$ column by column. Because the estimators of $\bw_i$, for $i=1,\cdots, p$, are obtained in the same manor, we will study the estimator for the $i$-th column and 
assume the number of non-zeros is exactly $s_0$ for simplicity.

Let $\bw_{i1}$ be the vector of non-zero elements of $\bw_i$, and then define $(\bB_i)_{p\times s_0}$ as a known matrix with either 0's or 1's such that $\bB_i \bw_{i1} = \bw_i$. 
In addition, we let $\bI_p$ be the $p \times p$ identity matrix,  $\bE_i \in \mathbb{R}^{p \times 1}$ is the $i$-th column of $\bI_p$. Let $\bI_{s_0}$ 
be the $s_0 \times s_0$ identity matrix, $\bg_i \in \mathbb{R}^{s_0 \times 1}$ is the $i$-th column of $\bI_{s_0}$.  
Our estimator is defined through the sample covariance matrix.  Let $\bar{\bX}=\sum_{k=1}^n \bX_k/n$ be the sample mean and $\bar{\bX} =(\bar{X}_1,\bar{X}_2,\cdots,\bar{X}_p)^\prime$. 
Let $\bSXn=\sum_{k=1}^n  (\bX_k - \bar{\bX})(\bX_k - \bar{\bX})^T/(n-1)$ be the sample covariance matrix. 

Let $\bsig_i$ denote the $i$-th sub-matrix of $\bsig$ corresponding to the non-zero elements of the $i$-th column $\bw_i$ of $\bOmg$ such that $\bB_i^T \bsig \bB_i$ = $\bsig_i$.
The inverse of $\bsig_i$ is then given by $(\bB_i^T \bsig \bB_i)^{-1} = \bsig_i^{-1} = \bOmg_i$. We will define the corresponding sample analogue of $\bsig_i$ as follows.
For $i = 1,\cdots, p$, let $\xbi = \bB_i^T \bX$ denote the sub-vector of $\bX$. 
Let $\bB_i^T \bX_1 ,\cdots, \bB_i^T \bX_n$ be IID copies of  $\bB_i^T \bX$, and denote them by $\x1bi,\cdots, \xnbi$, respectively.  
Note that $\xbi$ follows a multivariate normal distribution with mean vector $\bB_i^T\bmu \in \mathbb{R}^{s_0 \times 1}$ and covariance matrix $\bB_i^T \bsig \bB_i$. 
Let $\bar{\bX}_{\bB_i}  = \sum_{k=1}^n \xkbi/n$ be the sample mean of $\{\bX_{k\bB_i}\}_{i=1}^n$, and $\bSXni = \sum_{k=1}^n (\xkbi - \bar{\bX}_{\bB_i})(\xkbi - \bar{\bX}_{\bB_i})^T/(n-1)$ be the sample 
covariance matrix corresponding the non-zero components of the $i$-th column of $\bOmg$.

To obtain an estimator for $\bw_i$, we first represent $\bw_i$ using the population covariance matrix $\bsig$. Using the definition of $\bOmg$, $\bsig \bOmg = \bI_p$, we have, for $i=1,\cdots, p$, 
 \begin{equation}
 \label{equi} \bsig \bw_i = \bsig \bB_i\bw_{i1}  = \bE_i.
 \end{equation}

\begin{lemma}
\label{proestimator}
The unique solution $\bw_{i1}$ of equation (\ref{equi}) is $(\bB_i^T \bsig \bB_i)^{-1} \bB_i^T \bE_i$.
\end{lemma}

\noindent\textit {Proof}: Multiplying $\bB_i^T$ on the both sides of the equation (\ref{equi}) yields  $\bB_i^T \bsig \bB_i\bw_{i1}
= \bB_i^T \bE_i.$  Because $\bB_i^T \bsig \bB_i$ is an $s_0\times s_0$ invertible matrix, $\bw_{i1} = (\bB_i^T \bsig \bB_i)^{-1}\bB_i^T \bE_i$ 
 \hfill $\square$

Based on Lemma \ref{proestimator}, a natural plug-in estimator of $\bw_{i1}$ is $\hat{\bw}_{i1} = (\bB_i^T \bSXn \bB_i)^{-1} \bB_i^T \bE_i= \bSXni^{-1} \mathbf{f}_i$, where $i =1,\cdots, p.$ Accordingly, $\bw_i$ is estimated by 
$$\hat{\bw}_{i} = \bB_i(\bB_i^T \bSXn \bB_i)^{-1} \bB_i^T \bE_i =  \bB_i \bSXni^{-1}\mathbf{f}_i.$$
where $ \mathbf{f}_i = \bB_i^T \bE_i \in \mathbb{R}^{s_0 \times 1}$ and $ \mathbf{f}_i$ is a vector with one element being 1 and all the other elements being 0, and we
use the fact that $\bB_i^T \bSXn \bB_i$ = $\bSXni$.
Because the proposed estimator $\hat{\bw}_{i}$ is invariant to $\bmu$, without loss of generality, assume that $\bmu=0$ in the rest of the paper.

Beyond the class of Gaussian random variables, the proposed estimator can also be applied to other classes of distributions whenever a connection between the precision matrix and its graphical structure is available, for instance, the class of nonparanormal variables defined in Liu et al. (2009). 

\section{Main results}

In this section, we will study the rates of convergence and the asymptotic distribution of the proposed estimator. 
To establish consistency of the proposed estimators, we assume the following regularity conditions:

\begin{itemize}
\item[(C1)]  Let $0 < \lambda_1(\bsig) \leq \lambda_2(\bsig) \leq \cdots \leq \lambda_p(\bsig)$ 
be eigenvalues of $\bsig$. There exists some constant $C_0 >0$ such that $C_0 \leq \lambda_1(\bsig) \leq \lambda_p(\bsig) \leq C_0^{-1}$. 
\end{itemize}
\noindent \textit{Remark 1.}  Condition (C1) is a mild condition commonly used in the literature. For example, Zhou et al. (2011) and Liu and Luo (2015).

The following theorem gives us the proposed estimator's rate of convergence when we estimate nonzero positions of the precision matrix of one column.
Theorem \ref{columnrate} shows that the estimator $\hat{\bw}_{i1}$ for each column of precision matrix $\bOmg$ is consistent whence 
$s_0 = o(n)$.

 \begin{theorem}\label{columnrate} 
Under condition (C1), $ \v \hat{\bw}_{i1} - \bw_{i1} \v = {O_P}(\sqrt {{s_0}/{n}} )$, for $i =1, 2, \cdots, p$.
\end{theorem}

\noindent\textit{Proof:}  For each submatrix $\bsig_i$ of $\bsig$, suppose its eigenvalues are 
$0 \leq \lambda_1 (\bsig_i) \leq \cdots \leq \lambda_{s_0}(\bsig_i)$. Under condition (C1), we have 
$C_0 \leq \lambda_1 (\bsig_i) \leq \lambda_{s_0} (\bsig_i) \leq C_0^{-1} $, because of the eigenvalues interlacing theorem.
So the conditions in Lemma \ref{DisCov}-\ref{ine1} in the Appendix are satisfied and they can applied for every column estimator $\hat{\bw}_{i1}$.
Define $\bWXni = \bSXni - \bsig_i$.
Because $\v \hat{\bw}_{i1} - \bw_{i1} \v =\v  \bSXni^{-1}\mathbf{f}_i -  \bsig_i^{-1} \mathbf{f}_i \v$ and $\v \mathbf{f}_i\v=1$,
applying Lemma \ref{ine1} in the Appendix, we have
$$\v  \bSXni^{-1} \mathbf{f}_i -  \bsig_i^{-1} \mathbf{f}_i \v  \leq \v \mathbf{f}_i\v \v \bSXni^{-1} -  \bsig_i^{-1} \v \leq \v \bsig_i^{-1} \v  \frac{\v  \bsig_i^{-1}\bWXni \v}{1 - \v \bsig_i ^{-1} \bWXni \v}=O_P(\sqrt{{s_0}/{n}}),$$ 
where we applied Lemma \ref{rate1} in the Appendix, which implies that
$\v \bsig_i ^{-1} \bWXni \v = O_P(\sqrt{s_0/n}) \rightarrow 0$, as $n \rightarrow \infty$
because $s_0 = o(\sqrt{n}$). This implies that the condition $1 - \v \bsig_{i}^{-1} \bWXni \v >0$  in Lemma \ref{ine1} always holds asymptotically. 
Theorem \ref{columnrate} is proved \hfill $\square$

The entire precision matrix estimator of $\bOmg$ is $\hat{\bOmg}=(\hat{\bw}_1,\cdots, \hat{\bw}_p)$ where $\hat{\bw}_i$ can be 
obtained by $\hat{\bw}_i=\bB_i\hat{\bw}_{i1}$ for $i=1,\cdots, p$. Theorem \ref{uniformrate} provides the rate
of convergence of $\hat{\bOmg}$ in $L_1$ norm and the maximum of the $L_2$ norm. The proof of Theorem \ref{uniformrate} can be found in the Appendix.

\begin{theorem} 
\label{uniformrate} 

Under condition (C1), we have the uniform rate of convergence of the proposed estimators $\hat{\bw}_{i1}$ is 
$\max_{i=1,\cdots,p}\v \hat{\bw}_{i1} - \bw_{i1} \v = O_P(s_0\sqrt{{\log p}/{n}})$ and if $\v \bOmg \v_1 <\infty$ we have $\|\hat{\bOmg}-\bOmg\|_1=O_p(s_0\sqrt{{\log p}/{n}})$.

\end{theorem}

The above rate of convergence is the same as that established in Theorem 4.1 of Liu and Wang (2017). It is a minimax optimal rate (Liu and Luo, 2015) over precision matrices satisfying
condition (C1).  In the following theorem, we establish the rate of convergence of $\hat{\bOmg}$ in elementwise sup-norm rate norm, which is also minimax rate optimal (Liu and Luo, 2015).
 
\begin{theorem} \label{maxcomp}
Under conditions (C1) and if $\v \bOmg \v_1 <\infty$,  then $\max\limits_{i,j =1,\cdots,p} \l \hat{\omega}_{ij} - \omega_{ij} \l = O_P(\sqrt{{\log p}/{n}})$.
\end{theorem}

We now establishes the asymptotic normality of the linear combinations of the proposed estimators $\hat{\bw}_{i1}$ for $i=1 \cdots, p$. The following additional condition
is needed.

\begin{itemize}
\item[(C2)] There exists some constants $K_1 > 0$ such that $\v \bsig \v_1 \leq K_1$.
\end{itemize}

\noindent \textit{Remark 2.} Condition (C2) is commonly used in the literature. For example, conditions (C1) and (C2) are met by the class of covariance matrices as defined in page 6 of Bickle and Levina (2008).

 \begin{theorem} \label{sample.pre.asymp}
Let $\bm_{s_0 \times 1}=(m_1,m_2,\cdots, m_{s_0})^T$ and $\bm^T \bOmg_i = (a_1,\cdots, a_{s_0})$. Under conditions (C1) and (C2) and assume $a_k = O(1),$ for $k =1,2,\cdots,s_0$,  then
 $$\frac{\sqrt{n}}{\sqrt{h}}\bm^T (\hat{\bw}_{i1} - \bw_{i1}) \stackrel{d}{\to} N(0,1),$$
 where $h = \Var( \bm^T \bOmg_i \bX_{1i} \bX_{1i}^T \bOmg_i \mathbf{f}_i)$ for $i = 1,\cdots, p$.
 \end{theorem}
 
\noindent\textit{Proof:} From Lemma \ref{represent} in the Appendix, we have 
 \begin{align*}
\frac{\sqrt{n}}{\sqrt{h}}\bm^T (\hat{\bw}_{i1} - \bw_{i1})=\frac{\sqrt{n}}{\sqrt{h}}\bm^T (\bSXni^{-1} - \bOmg_i)\mathbf{f}_i  
&= -\frac{\sqrt{n}}{\sqrt{h}} \bm^T \bOmg_i \bWXni \bOmg_i \mathbf{f}_i 
 -\frac{\sqrt{n}}{\sqrt{h}} \bm^T \bOmg_i \bWXni (\bSXni^{-1}- \bOmg_i) \mathbf{f}_i \\
 &= I_1 - I_2.
 \end{align*}

Suppose that $\bOmg_i \mathbf{f}_i  =  (b_1,\cdots, b_{s_0})^T$. 
We have $b_k = O(1)$ for $k = 1,\cdots, s_0$. In addition, 
 $a_k = O(1).$ So $a_ib_j = O(1)$ for any $1 \leq i, j \leq s_0$. 
 Applying Lemma \ref{sample.cov.asymp} in the Appendix to $I_1$, we have
 $$-\frac{\sqrt{n}}{\sqrt{h}} \bm^T \bOmg_i \bWXni \bOmg_i \mathbf{f}_i
 = -(\frac{\sqrt{n}}{\sqrt{h}} \ba^T \bSXni \bb - \frac{\sqrt{n}}{\sqrt{h}} \ba^T \bsig_i \bb) \stackrel{d}{\to} N(0,1).$$
The theorem is proved if we can show that $I_2$ is ignorable. 

Lemma \ref{resrate} in the Appendix implies that
 $ \l \bm^T \bOmg_i \bWXni (\bSXni^{-1}- \bOmg_i) \mathbf{f}_i \l 
 = O_P({s_0 \v \bm \v}/{n})$.
Lemma \ref{variance} in the Appendix gives us
 $h  \asymp O(\bm^T \bOmg_i \bm)$.  
Combining these two together gives us
 $$\frac{\sqrt{n}}{\sqrt{h}} \l \bm^T \bOmg_i \bWXni (\bSXni^{-1}- \bOmg_i) \mathbf{f}_i \l 
 = O_P(\frac{s_0 \v \bm \v}{\sqrt{n} \sqrt{\bm^T \bOmg_i \bm} }).$$
 
 Suppose $0 < \lambda_1(\bOmg_i) \leq \cdots \leq \lambda_{s_0}(\bOmg_i)$ are eigenvalues of $\bOmg_i$, by eigenvalues interlacing theorem, we have  
 $C \leq \lambda_1(\bOmg_i) \leq \lambda_{s_0}(\bOmg_i) \leq C^{-1}, $ for some constant $C >0$, and
$ i =1,\cdots,p$. Notice that $C \v  \bm \v^2 = \epsilon_0 \v \bm^T \bm \v < \lambda_1(\bOmg_i) \v \bm^T \bm \v  \leq \bm^T \bOmg_i \bm  $
  and $\bm^T \bOmg_i \bm \leq \lambda_{s_0}(\bOmg_i) \v \bm^T \bm \v \leq C^{-1} \v \bm^T \bm \v= C^{-1} \v \bm \v^2.$
  So  $\sqrt{{\bm}^T {\bOmg}_i {\bm}} \asymp \v {\bm} \v$. This gives us $\sqrt{n}\l {\bm}^T {\bOmg}_i {\bWXni} (\bSXni^{-1}- {\bOmg}_i) {\bf}_i \l/\sqrt{h}= O_P({s_0}/{\sqrt{n}} ) = O_P(1).$
  In other words, $I_2$ is ignorable. The theorem is proved \hfill $\square$

%
%
%

 \section{Simulation Studies}
In this section, we perform some numerical analysis to illustrate the proposed estimator's finite sample performance. We conduct two sets of simulation studies. In the first set of simulations, the goal is to demonstrate the proposed estimator's accuracy. In the second set of simulation studies, we compare the proposed method to another method proposed by Liu and Wang (2017), 
which does not utilize the network structure. Because the proposed method assumes some prior knowledge of the graphical structure, the method of Liu and Wang (2017) is used as a benchmark to check
the performance of the proposed estimator.

We generate $n$ independent and identically distributed multivariate normal distributed $p$-dimensional random vectors $N(0,\bsig)$ with mean vector 0 and covariance matrix $\bsig$, where $\bsig^{-1} =\bOmg=(\omega_{ij})_{p \times p}$  
such that $\omega_{ij}= 0.6^{\l i-j \l}$ for  $\l i-j \l < s_0$ and $\omega_{ij}=0$ other wise. 
Because we estimate the precision matrix $\bOmg$ in column by column fashion, for simplicity, we only demonstrate the accuracy of our estimator for the first column of the precision matrix. 
Note that the first column of $\bOmg$ has $s_0$ non-zero components $\bw_{11} =(\omega_{11,1}, \omega_{11,2}, \cdots, \omega_{11,s_0})^T=(0.6^0, 0.6^1,\cdots, 0.6^{s_0-1})^T.$

\subsection{Accuracy of the Proposed Estimator}

To evaluate the proposed estimator's accuracy, we report the bias of the proposed estimators for the non-zero components $\omega_{11,k}$ ($k = 1,2,\cdots, s_0$) and standard deviations of the bias. 
Both relative bias and absolute bias are evaluated in the simulation. To summarize the overall performance of the estimation of the entire non-zero vector $\omega_{11}$, we define the relative bias
and the absolute bias, respectively, as follows:
$$\mbox{RelBias} =\frac{1}{s_0}\sum\limits_{k=1}^{s_0} \frac{\l \hat{\omega}_{11,k} - \omega_{11,k} \l}{|\omega_{11,k}|},\;\;\mbox{and}\;\; \mbox{AbsBias} =\frac{1}{s_0}\sum\limits_{k=1}^{s_0} \l \hat{\omega}_{11,k}-\omega_{11,k}\l. $$
In addition to the overall performance, we also report the bias and its standard deviation of the component-wise estimators. All the simulation results reported in this section are based on 300 simulation replications. 
To understand the effect of sample size and data dimension,  we consider three different sample sizes $n = 100, 300,$ and 500. For each sample size $n$, we change the data dimension according to
the ratios ${p}/{n}$ at five different values $0.1, 0.5, 1, 5,$ and 10. 

Table \ref{tab1:title} reports the average and standard deviations of the RelBias and AbsBias, and the average and standard deviations of the non-zero components in $\bw_1$ with sparsity level at $s_0=4$. 
Note that in this simulation, $\bw_{11} = (\omega_{11,1}, \omega_{11,2},\omega_{11,3}, \omega_{11,4})^T=(1, 0.6,0.36, 0.216)^T.$
The standard deviations are included in the parentheses. We observe that the data dimension has no significant impact on the proposed estimator. As the sample size increases from 100 to
500, the relative bias and absolute bias of the proposed estimator decrease significantly. This illustrates the consistency of the proposed estimator.  Similar phenomena can be observed from
the component-wise estimators for $\omega_{11,1}, \omega_{11,2},\omega_{11,3}$, and $\omega_{11,4}.$
\begin{table}
\begin{center}
\captionof{table}{Relative bias, absolute bias and the standard deviations of the proposed estimator for the first column $\bw_{11}$ of $\bOmg$. 
The average and standard deviation of the estimators for the components in $\bw_{11}$ are also included. The standard deviations are included in
the parentheses.} \label{tab1:title} 
\begin{tabular}{crrrrrrr}
  \hline
$n$ & ${p}/{n}$ & RelBias & AbsBias & $\omega_{11,1}$ & $\omega_{11,2}$& $\omega_{11,3}$ & $\omega_{11,4}$\\
  \hline
100 &   0.1 & 0.21 (0.12) & 0.09 (0.06) & 1.05 (0.15) & 0.62 (0.12) & 0.38 (0.11) & 0.23 (0.08) \\ 
 &   0.5 & 0.22 (0.13) & 0.10 (0.06) & 1.05 (0.16) & 0.63 (0.13) & 0.38 (0.11) & 0.23 (0.09) \\ 
 &   1.0 & 0.20 (0.12) & 0.09 (0.06) & 1.05 (0.15) & 0.63 (0.13) & 0.38 (0.11) & 0.22 (0.08) \\ 
 &   5.0 & 0.21 (0.12) & 0.09 (0.05) & 1.05 (0.14) & 0.64 (0.12) & 0.38 (0.11) & 0.23 (0.09) \\ 
 & 10.0 & 0.22 (0.13) & 0.10 (0.06) & 1.06 (0.16) & 0.63 (0.13) & 0.38 (0.11) & 0.23 (0.09) \\ 
   \hline
300  &  0.1 & 0.11 (0.06) & 0.05 (0.03) & 1.02 (0.08) & 0.61 (0.06)   & 0.37 (0.06) &  0.22 (0.05) \\ 
  &  0.5   & 0.11 (0.06) & 0.05 (0.03) & 1.01 (0.08) & 0.61 (0.06) & 0.36 (0.06) &  0.22 (0.05) \\ 
  &  1.0   & 0.12 (0.06) & 0.05 (0.03) & 1.02 (0.08) & 0.61 (0.07)    & 0.36 (0.06) &  0.22 (0.05) \\ 
  &  5.0   & 0.12 (0.06) & 0.05 (0.03) & 1.02 (0.09) & 0.61 (0.07)    & 0.36 (0.06) &  0.22 (0.05) \\ 
  &  10.0 & 0.12 (0.06) & 0.05 (0.03) & 1.02 (0.09) & 0.62 (0.07)   & 0.37 (0.05) &  0.23 (0.05) \\ 
   \hline
500  &  0.1 & 0.09 (0.05) & 0.04 (0.02) & 1.02 (0.07) & 0.61 (0.05)  &  0.37 (0.05)  & 0.22 (0.04) \\ 
  &  0.5    & 0.09 (0.05) & 0.04 (0.02) & 1.01 (0.06) & 0.61 (0.05)  &  0.37 (0.05)  & 0.22 (0.04) \\ 
  &     1.0 & 0.09 (0.04) & 0.04 (0.02) & 1.01 (0.07) & 0.61 (0.05)  &  0.36 (0.04)  & 0.22 (0.03) \\ 
  &     5.0 & 0.09 (0.05) & 0.04 (0.02) & 1.01 (0.06) & 0.61 (0.05)  &  0.36 (0.04)  & 0.22 (0.04) \\ 
  &   10.0 & 0.09 (0.05) & 0.04 (0.02) & 1.02 (0.07) & 0.61 (0.05)  &  0.37 (0.05)  & 0.22 (0.03) \\ 
   \hline
\end{tabular}
\end{center}
\end{table}

To illustrate the impact of the sparsity level on the proposed estimator, we increased the sparsity level $s_0$ from 4 to 6 in the following Table \ref{tab2:title}. In this scenario, we need to estimate 6 coefficients in the first 
column of the precision matrix, which is $\bw_{11}  = (\omega_{11,1}, \omega_{11,2}, \omega_{11,3}, \omega_{11,4}, \omega_{11,5}, \omega_{11,6} )^T = (1, 0.6, 0.36, 0.216, 0.1296, 0.07776)^T$. 
Table \ref{tab2:title} shows that the increasing of sparsity level $s_0$ has some effect on the accuracy of the proposed estimator. As the number of non-zero components increases, the accuracy of the proposed estimator
drops. This is understandable due to the increasing in the number of unknown parameters. The other patterns observed in Table \ref{tab2:title} are very similar to those in Table \ref{tab1:title}. 

\begin{table}
\tabcolsep 4pt
\centering
\captionof{table}{ Relative bias, absolute bias and the standard deviations of the proposed estimator for the first column $\bw_{11}$ of $\bOmg$. 
The average and standard deviation of the estimators for the components in $\bw_{11}$ are also included. The standard deviations are included 
in the second row in every scenario.} \label{tab2:title} 
\begin{tabular}{crrrrrrrrr}
  \hline
n&${p}/{n}$ & RelBias  & AbsBias  & $\omega_{11,1}$  & $\omega_{11,2}$ & $\omega_{11,3}$ & $\omega_{11,4}$ & $\omega_{11,5}$ & $\omega_{11,6}$  \\  
  \hline
100 &  0.1   & 0.41  & 0.1  & 1.07  & 0.64  & 0.38 & 0.23  & 0.14  & 0.08  \\
       &          & (0.19) & (0.05) & (0.16) & (0.14) & (0.12) & (0.1) & (0.1) & (0.09) \\   
       &   0.5  & 0.45  & 0.1  & 1.07  & 0.65  &  0.4  &  0.24  & 0.15 & 0.08  \\ 
       &          &  (0.22) & (0.05) & (0.16) & (0.13) &  (0.11) &  (0.11) & (0.11) & (0.1)  \\
       &      1  & 0.41  & 0.1  & 1.08  & 0.64  &  0.38   & 0.22  & 0.13  & 0.08 \\ 
       &          & (0.2) & (0.05) & (0.16) & (0.13) &  (0.12) &   (0.1) & (0.1) & (0.1) \\ 
       &      5  & 0.42  & 0.09  & 1.07  & 0.64  & 0.39 &  0.23  & 0.15 & 0.09 \\
       &          &  (0.2) & (0.04) & (0.15) & (0.12) & (0.11) &  (0.11) &  (0.1) & (0.09)\\
       &     10 & 0.44  & 0.1  & 1.08  & 0.65  & 0.39 & 0.23  & 0.14  & 0.08  \\ 
       &          &  (0.21) & (0.05) & (0.16) & (0.14) & (0.12) & (0.11) & (0.1) & (0.1)\\ 
   \hline
300 &  0.1 & 0.23 & 0.05  & 1.0 & 0.61  & 0.37 & 0.22  & 0.13  & 0.08   \\
       &         & (0.1) & (0.02) & (0.08) & (0.07) & (0.06) & (0.06) & (0.06) & (0.05)  \\ 
       &  0.5 & 0.22 & 0.05  & 1.02 & 0.61  & 0.37 &  0.22  & 0.14  & 0.08  \\
       &        &  (0.1) & (0.02) & (0.09) & (0.07) & (0.06) &  (0.06) & (0.06) & (0.05)  \\
       &     1 & 0.22  & 0.05  & 1.02  & 0.61 & 0.36 & 0.21  & 0.13  & 0.08   \\ 
       &        & (0.1) & (0.02) & (0.08) & (0.07)& (0.06)& (0.06) & (0.06) & (0.05) \\ 
       &      5 & 0.22  & 0.05  & 1.02  & 0.61 & 0.36 & 0.22 & 0.13 & 0.08  \\
       &        &  (0.11) & (0.02) & (0.08) & (0.07) & (0.06) & (0.06) & (0.06) & (0.05) \\ 
       &    10 & 0.23 & 0.05 & 1.02 & 0.62  &  0.37 & 0.22  & 0.13  & 0.08  \\ 
       &        &  (0.11) & (0.02) & (0.09) & (0.07) &  (0.07) & (0.06) & (0.06) & (0.05)\\ 
   \hline
500 &  0.1 & 0.18  & 0.04 & 1.02  & 0.61   & 0.36   &  0.22  & 0.13  & 0.08  \\ 
       &         & (0.09) & (0.02) & (0.06) & (0.05)  & (0.05)  &  (0.05) & (0.04) & (0.04) \\ 
       &  0.5 & 0.17  & 0.04 & 1.02  & 0.61   & 0.37   &  0.22  & 0.13  & 0.08  \\ 
       &        & (0.08) & (0.02) & (0.06) & (0.05)  & (0.05)  &  (0.05) & (0.04) & (0.04) \\ 
       &     1 & 0.18  & 0.04 & 1.01  & 0.61   & 0.37   &  0.23  & 0.14  & 0.09  \\ 
       &        & (0.09) & (0.02) & (0.06) & (0.05)  & (0.05)  &  (0.05) & (0.04) & (0.04) \\ 
       &     5 & 0.18  & 0.04 & 1.01  & 0.61   & 0.36   &  0.22  & 0.13  & 0.08  \\ 
       &        & (0.09) & (0.02) & (0.07) & (0.05)  & (0.05)  &  (0.04) & (0.05) & (0.04) \\ 
       &   10 & 0.17  & 0.04 & 1.02  & 0.61   & 0.37   &  0.22  & 0.13  & 0.08  \\ 
       &        & (0.07) & (0.02) & (0.07) & (0.05)  & (0.05)  &  (0.04) & (0.04) & (0.03) \\ 
   \hline
\end{tabular}
\end{table}

\subsection{Comparison with Some Existing Estimators}

This subsection compares the proposed estimator with the TIGER estimator proposed by Liu and Wang (2017).  The proposed estimator assumes that the graphical structure is known and incorporates the information in formulating the estimator, but the TIGER estimator does not use any knowledge of the graphical structure. Our goal is to check the improvement of utilizing the known graphical structure to construct the precision matrix estimator.

We briefly introduce the TIGER approach to estimate the precision matrix. TIGER stands for Tuning-Insensitive Graph Estimation and Regression, proposed by Liu and Wang (2017). This estimating procedure for the precision matrix uses a column by column approach. TIGER uses the similar idea as Meinshausen and Buhlman (2006), which estimates the unknown precision matrix by regressing $i$-th variable ($i = 1,\cdots,p$)
on the remaining $p-1$ variables. However, unlike Meinshausen and Buhlman (2006) approach, each subproblem is solved by using the square root-LASSO (SQRT-LASSO). To be more explicit, let
$\bZ_{n \times p} =(Z_{ij})_{i=1,j=1}^{n,p}$
be the normalized data matrix where $n$ is the sample size, $p$ is the dimension, and the sample variance of each variable is one. Let $\bZ_{i*} \in \mathbb{R}^{1 \times n}$ is the $i$-th row of $\bZ$, $\bZ_{*j} \in \mathbb{R}^{p \times 1}$ is the $j$-th column of $\bZ$, $\bZ_{\mbox{\textbackslash} j,j} \in \mathbb{R}^{(p-1) \times 1}$ corresponds to the $j$-th column of $\bZ$ with the $j$-th row removed, $\bZ_{* \mbox{\textbackslash}j} \in \mathbb{R}^{p \times (n-1)}$ is the submatrix of $\bZ$ with the $j$-th column removed, and $\bZ_{\mbox{\textbackslash}j, \mbox{\textbackslash}j} \in \mathbb{R}^{(p-1) \times (n-1)}$ is the submatrix of $\bZ$ with the $j$-th row and the $j$-th column removed.
Then for $j = 1,\cdots,p$, $j$-th column of the precision matrix $\bOmg$ is estimated by $\hat{\bOmg}_{jj} = \hat{\tau}_j^{-2} \hat{\bgam}_{jj}^{-1}$ and $\hat{\bOmg}_{\mbox{\textbackslash} j,j} =  - \hat{\tau}_j^{-2} \hat{\bgam}_{jj}^{-1/2}\hat{\bgam}_{\mbox{\textbackslash} j,\mbox{\textbackslash} j}^{-1/2}\hat{\bbeta}_j$, where $\hat{\bbeta}_j =  \argmin_{\bbeta_j \in R^{p-1}} \{ \v \bZ_{*j} - \bZ_{* \mbox{\textbackslash}j} \bbeta_{j} \v_2/\sqrt{n} + \lambda \v \bbeta_j \v_1 \}$, $\hat{\tau}_j = \v \bZ_{*j} - \bZ_{* \mbox{\textbackslash}j} \hat{\bbeta}_{j} \v_2/\sqrt{n}$, and $\bgam = \text{diag}(\bS_{\bZ,n})$ is the diagonal matrix with elements are the main diagonal elements of its sample covariance matrix. As we can see here, when estimating parameter $\bbeta_{j}$, TIGER uses SQRT-LASSO instead of LASSO penalty as in Meinshausen and Buhlman (2006). When the data matrix is not normalized, we can refer to Liu and Wang(2017) for more details. The TIGER method's main advantage compared to the other methods is that it is asymptotic free of tuning and does not depend on any unknown parameters due to the SQRT-LASSO property (Belloni et al., 2012).  Nevertheless, in reality, the sample size is always finite. So an optimal estimation for precision matrix for a finite sample is still chosen by comparing loss functions at different tuning values via some selection tuning values approaches, such as K fold cross validation or StARS (Stability Approach to Regularization Selection). In the below simulation, we used the ``flare'' R package created by the authors to obtain the optimal precision matrix. We applied the K fold cross-validation approach to obtain the optimal precision matrix under the default setting in the package. In this default setting tuning parameters is a sequence of 5 values uniformly distributed in the interval $[\pi \sqrt{\log p /n}/4, \pi \sqrt{\log p /n}]$. The process of choosing the optimal precision matrix is processed as follows. At each tuning value $\lambda$, the precision matrix $\hat{\bOmg}_\lambda$ is estimated from the training sets, the likelihood loss function $tr(\bS_{\bZ,n})\hat{\bOmg}_\lambda - 
\log\{\det(\hat{\bOmg}_\lambda)\}$ is calculated based on the leave out data fold. Then the optimal precision matrix is corresponding to the $\lambda$ which yields the smallest value for average absolute loss. Now we are ready to compare the performance of our proposed method and the benchmark method TIGER.

Data are generated according to the same setting described in Subsection 4.1. Similar to the previous subsection, we consider two sample sizes $n = 100$ and 200, five dimension and sample size ratios: 
$p/n = 0.1, 0.5, 1, 5,$ and $10$, and two sparsity levels $s_0 = 4$ and $s_0 = 6$.  The comparison is made through comparing the relative bias and absolute bias defined in Subsection 4.1.

\begin{table}
\begin{center}
\captionof{table}{Relative bias, absolute bias and their standard deviations of the proposed estimator and the TIGER method for the first column $\bw_{11}$ of $\bOmg$. 
The standard deviations are all included in the parentheses. \label{tigervsprop}}
\begin{tabular}{ccrcccc}
\hline
     &    &  &  \multicolumn{2}{c}{Proposed Estimator} & \multicolumn{2}{c}{TIGER Estimator}\\
$s_0$ & $n$ & ${p}/{n}$  &  RelBias & AbsBias & RelBias  & AbsBias \\ 
  \hline    
4 & 100 &    0.1 &  0.24 (0.13) & 0.11 (0.06) &  0.65 (0.06) & 0.25 (0.04) \\
 &        &    0.5 &  0.21 (0.12) & 0.09 (0.05) &  0.69 (0.05) & 0.28 (0.04) \\
 &        &       1 &    0.2 (0.11) & 0.09 (0.05) &  0.71 (0.05) & 0.29 (0.04) \\
 &        &       5 &  0.23 (0.15) &   0.1 (0.07) &  0.73 (0.05) & 0.30 (0.04) \\
 &        &     10 &    0.2 (0.09) & 0.09 (0.04) &  0.74 (0.05) & 0.31 (0.04) \\
 \cline{2-7}
 & 200 &    0.1 &  0.14 (0.07) & 0.07 (0.03) &   0.62 ( 0.06) & 0.24 (0.03)  \\
 &        &    0.5 &  0.15 (0.1) & 0.07 (0.05) &  0.66( 0.04) & 0.26 (0.03)  \\
 &        &       1 &  0.14 (0.07) & 0.06 (0.03) &   0.68 (0.04) & 0.27 (0.03) \\
 &        &       5 &  0.16 (0.08) & 0.07 (0.04) & 0.70 (0.04) & 0.28 (0.03)  \\
 &        &     10 &  0.15 (0.08) & 0.07 (0.04) &  0.71 (0.04) & 0.29 (0.03) \\
 \hline   
6  & 100 &  0.1 &  0.45 (0.26) & 0.1 (0.06)  &   0.77 (0.04) & 0.2 (0.02)  \\
 &        &  0.5 &  0.44 (0.22) & 0.1 (0.05)  &   0.8 (0.03) & 0.22 (0.02)  \\
 &        &     1 &      0.43 (0.2) & 0.1 (0.05) &    0.81 (0.03) & 0.23 (0.03) \\
 &        &     5 &    0.4 (0.19) & 0.09 (0.05) & 0.82 (0.03) & 0.23 (0.02)  \\
 &        &   10 &   0.46 (0.24) & 0.1 (0.05) &   0.82 (0.03) & 0.24 (0.02) \\
 \cline{2-7}
 & 200 &  0.1 &   0.3 (0.14) & 0.07 (0.03) & 0.75 (0.04) & 0.19 (0.02)\\ 
 &        &  0.5 &  0.31 (0.17) & 0.07 (0.04) &  0.77 (0.03) & 0.21 (0.02)\\ 
 &        &     1 &  0.31 (0.17) & 0.07 (0.03)&   0.78 (0.03) & 0.21 (0.02) \\ 
 &        &     5 &   0.28 (0.13) & 0.06 (0.02) & 0.8 (0.03) & 0.22 (0.02)  \\ 
 &        &   10 &   0.28 (0.13) & 0.06 (0.03) &   0.8 (0.02) & 0.22 (0.02) \\
 \hline 
\end{tabular}
\end{center}
\end{table}

Table \ref{tigervsprop} summarizes the performance of the proposed estimator and the TIGER estimator. 
To make an informative comparison, the relative bias and absolute bias are only compared for the non-zero components in $\bw_1$. 
Based on Table \ref{tigervsprop}, we observe that the relative bias and absolute bias of the TIGER approach are much higher than the proposed method. This observation implies that
the proposed estimator, using the graphical structure information, outperforms the TIGER approach.  This result is something expected but confirmed with simulation studies. 
The results in Table \ref{tigervsprop} confirm that utilizing the graphical structure in estimating the precision matrix improves the estimator's accuracy.  

\section{Conclusion} 
  
We introduce a simple estimator for precision matrix when its graphical structure is known under the high-dimensional framework. The proposed estimator comes in an explicit form and requires no additional iterative 
algorithm. Also, the estimating procedure is constructed by a simple argument without using any likelihood function, which makes the method adapt to a non-parametric framework. The rates of convergence are obtained
under the $L_1$ norm and element-wise maximum norm where both rates are mini-max optimal. We also derive the asymptotic normality for a linear combination of the estimators for precision matrix columns. 
Asymptotic normality property makes statistical inference available with the proposed estimator. Simulation studies demonstrate that incorporating network structure in estimating precision matrix significantly improve the 
accuracy of the estimator.    

\section*{Acknowledgements}
This research was partially supported by the NSF grant DMS-1462156.

\appendix
\newpage
\noindent \textbf{APPENDIX}

Define $\bVXn=\sum_{k=1}^n \bX_k \bX_k^T/n$. The corresponding submatrix for the non-zero components of the $i$-th column of $\bOmg$ is $\bVXni =\sum_{k=1}^n \xkbi \xkbi^T/n$ 
such that $\bB_i^T \bVXn \bB_i = \bVXni$. In addition, we denote $\bWXni = \bSXni - \bsig_i$. 
Moreover, the components of matrices $\bsig$, $\bOmg $, $\bSXn$, $\bVXn$ and our estimator precision matrix $\hat{\bOmg}$ 
are represented as the following $\bsig = (\sigma_{ij})_{p \times p}$,  $\bOmg = (\omega_{ij})_{p \times p}$, $\bSXn = (s_{ij})_{p \times p}$, $\bVXn = (v_{ij})_{p \times p}$, and $\hat{\bOmg} = (\hat{\omega}_{ij})_{p \times p}$. We denote the components of  $\xkbi, \xbi$ and $\bar{\bX}_{\bB_i}$, respectively, by $\xkbi = (X_{k\bB_i, 1},\cdots, X_{k\bB_i, s_0})$ for $k =1,2,\cdots,n$ 
, $\xbi = (X_{\bB_i,1},\cdots, X_{\bB_i, s_0})$ and $\bar{\bX}_{\bB_i} =(\bar{X}_{\bB_i,1},\bar{X}_{\bB_i,2},\cdots,\bar{X}_{\bB_i,s_0})$.

To simplify our notations, we sometimes use a generic notation $\bY$ to denote the random vectors $\bX_{\bB_i}$, so our results can be applied to estimators of $\bw_i$ in each column. 
Let $\bY_1,\cdots, \bY_n$ be $s_0 \times 1$ random vectors, which are IID copies of $\bY = (Y_1, Y_2, \cdots, Y_{s_0})'$ from a multivariate normal distribution with mean vector 0 and covariance matrix $\bsig_{\bY}$
where $\bsig_{\bY} \in \mathbb{R}^{s_0 \times s_0}$, $\bsig_{\bY} > 0$ and $\bOmg_{\bY} = \bsig_{\bY}^{-1}$. Let $\bar{\bY}=\sum_{k=1}^n \bY_k/n$ be the sample mean vector, 
$\bar{\bY} =(\bar{Y}_1,\bar{Y}_2,\cdots,\bar{Y}_{s_0})'$. In addition, let $\bSYn=\sum_{k=1}^n (\bY_k-\bar{\bY})(\bY_k-\bar{\bY})^T/(n-1)$ 
be the sample covariance matrix, $\bVYn = \sum_{k=1}^n \bY_k \bY_k^T/n$, and $\bWYn = \bSYn - \bsig_{\bY}$.

Lemma \ref{DisCov} is a result on the rate of consistency of a sample covariance matrix, which was proved in Vershynyn (2011). 
The result in Lemma \ref{DisCov} holds for a larger class of distributions beyond the sub-Gaussian distribution, see Theorem 1 in Adamczak (2011) for more details.

\begin{lemma}
\label{DisCov}
Let $\lambda_1(\bsig_{\bY}) \leq \lambda_2(\bsig_{\bY}) \leq \cdots \leq \lambda_{s_0}(\bsig_{\bY})$ be eigenvalues of $\bsig_{\bY}$. If $\lambda_{s_0}(\bsig_{\bY})$ is bounded from the above, then 
$\v \bVYn - \bsig_{\bY} \v = O_P( \sqrt{{s_0}/{n}}).$
\end{lemma}

\noindent\textit{Proof:} Since Gaussian variables is a special case of sub-Gaussian, this is a direct result from Theorem 5.39 in Vershynyn (2011) \hfill $\square$ 

Note that  $\bSYn = n\bVYn/(n-1) + n\bar{\bY}\bar{\bY}^T/(n-1)$. Theorem 1 in Greif (2006) and Lemma \ref{DisCov} imply that
\begin{align}
   \v {\bSYn - \bsig_{\bY} } \v &\leq \v \frac{n}{n-1}\bVYn -\bsig_{\bY}  \v \leq \v \bVYn - \bsig_{\bY}  \v  + \v \frac{1}{n-1} \bVYn \v = O_P ( {\sqrt {{s_0}/{n}} } ).\label{consistent3}
\end{align}
In other words, (\ref{consistent3}) implies $\v \bWYn \v = O_P( \sqrt{{s_0}/{n}} )$ \hfill $\square$

\begin{lemma}
\label{rate1}
Let $0 < \lambda_{1}(\bsig_{\bY}) \leq \cdots \leq \lambda_{s_0}(\bsig_{\bY}) $ be eigenvalues of $\bsig_{\bY}$. If there are some constants $M, N >0$ such that 
$ N \leq \lambda_{1}(\bsig_{\bY}) \leq \lambda_{s_0}(\bsig_{\bY}) \leq M $. 
Then, we have $\v \bsig_{\bY}^{-1} \v = O(1)$ and $\v \bsig_{\bY}^{-1} \bWYn \v = O_P(\sqrt{{s_0}/{n}}).$
 \end{lemma}
 
 \noindent\textit{Proof:}
  Since $\lambda_{1}(\bsig_{\bY}) \leq \cdots \leq \lambda_{s_0}(\bsig_{\bY}) $ 
 are eigenvalues of $\bsig_{\bY}$, eigenvalues of $\bsig_{\bY}^{-1}$ are ${1}/{\lambda_1(\bsig_{\bY})}, {1}/{\lambda_2(\bsig_{\bY})},\cdots, {1}/{\lambda_{s_0}(\bsig_{\bY})}. $
Theorem 2.3.1 in Golub et al. (2015) tells us the operator norm of a symmetric matrix is the maximum eigenvalue of the matrix.
So $ \v \bsig_{\bY}^{-1} \v  = {1}/{\lambda _{1}(\bsig_{\bY})} \leq {1}/{N}$ or $\v \bsig_{\bY}^{-1} \v = O(1)$. 
We have $\v \bsig_{\bY}^{-1}\bWYn \v \leq \v \bsig_{\bY} ^{-1} \v \ \v \bWYn \v $ and $\v \bWYn \v= {O_P} \big ( {\sqrt {{s_0}/{n}} } \big )$. 
The lemma is proved \hfill $\square$

Under conditions of Lemma \ref{rate1},  we have $\v \bsig_{\bY}^{-1}\bWYn \v = O_P( \sqrt{{s_0}/{n}} ) = o_P(1)$. So, the condition $1 - \v \bsig_{\bY}^{-1} \bWYn \v >0$  in Lemma \ref{ine1} always holds 
 asymptotically. Then, a direct application of Lemma \ref{ine1} gives, for some $C>0$
 \begin{equation}\label{DisInv}
 \v \bSYn^{-1} -  \bsig_{\bY}^{-1} \v \leq C \v \bWYn \v = C \v \bSYn - \bsig_{\bY} \v.
 \end{equation}

\begin{lemma} \label{ine1} The following inequality about the sample covariance matrix inverse holds.
$$\v  \bSYn^{-1} -  \bsig_{\bY}^{-1} \v \leq \v \bsig_{\bY} ^{-1} \v \frac{\v  \bsig_{\bY}^{-1} \bWYn \v}{1 - \v \bsig_{\bY}^{-1} \bWYn \v},$$
 whenever  $1 - \v  \bsig_{\bY}^{-1} \bWYn \v > 0$.
\end{lemma}

\noindent \textit{Proof:} We have
  \begin{align}
  \v  \bSYn^{-1} \bWYn \v 
  &= \v \bSYn^{-1} \bsig_{\bY} \bsig_{\bY}^{-1} \bWYn \v =  \v \bSYn^{-1} \bsig_{\bY} \bsig_{\bY}^{-1} \bWYn 
  -  \bsig_{\bY}^{-1} \bWYn +  \bsig_{\bY}^{-1} \bWYn \v \notag\\
   &\leq  \v \bSYn^{-1}\bsig_{\bY} \bsig_{\bY}^{-1} \bWYn
    - \bsig_{\bY}^{-1}\bWYn \v + \v \bsig_{\bY}^{-1} \bWYn \v \notag
    \leq  \v  \bSYn^{-1} \bsig_{\bY}  - \bI_{s_0} \v  \v \bsig_{\bY}^{-1} \bWYn \v  +  \v  \bsig_{\bY}^{-1} \bWYn \v \notag.
  \end{align}
  
 Notice that $ \v \bSYn^{-1} \bWYn \v 
 = \v \bSYn^{-1} (\bSYn - \bsig_{\bY} ) \v 
 = \v \bI_{s_0} -  \bSYn^{-1} \bsig_{\bY}  \v. $
This yields $\v \bI_{s_0} - \bSYn^{-1} \bsig_{\bY}  \v\leq \v  \bSYn^{ -1} \bsig_{\bY}  - \bI_{s_0} \v \v \bsig_{\bY}^{-1} \bWYn \v + \v \bsig_{\bY}^{-1} \bWYn \v.$
It is equivalent to  
$\v \bI_{s_0} - \bSYn^{-1} \bsig_{\bY} \v  
\big ( 1 - \v \bsig_{\bY}^{-1} \bWYn \v \big )
 \leq  \v  \bsig_{\bY}^{-1} \bWYn \v,$
or $ \v \bI_{s_0} -  \bSYn^{-1} \bsig_{\bY}  \v\leq {\v  \bsig_{\bY} ^{-1} \bWYn \v}/\{1 - \v \bsig_{\bY}^{-1} \bWYn \v\},$
if $1 - \v \bsig_{\bY}^{-1} \bWYn \v>0 $.
It follows that
 $\v \bI_{s_0} -  \bSYn^{-1} \bsig_{\bY}  \v  \v \bsig_{\bY}^{-1} \v
  \leq \v \bsig_{\bY}^{-1} \v {\v  \bsig_{\bY}^{-1} \bWYn \v}/\{1 - \v \bsig_{\bY}^{-1} \bWYn \v\}.$
Then $\v \bsig_{\bY}^{-1} -  \bSYn^{-1} \v \leq \v \bsig_{\bY}^{-1} \v {\v  \bsig_{\bY}^{-1} \bWYn \v}/\{1 - \v \bsig_{\bY}^{-1} \bWYn \v\},$
thereby proving the lemma \hfill $\square$
 
\smallskip


\noindent \textit{Proof of Theorem 2:} We will only show the result $\max\limits_{i=1,\cdots,p}\v \hat{\bw}_{i1} - \bw_{i1} \v = O_P(s_0\sqrt{{\log p}/{n}})$.
since the proof for $\|\hat{\bOmg}-\bOmg\|_{1}=\max\limits_{i=1,\cdots,p} \sum_{j=1}|\hat{\bw}_{ij} - \bw_{ij}|= O_P(s_0\sqrt{{\log p}/{n}})$ is very similar.
For any given $t$, we have
\begin{align}
\label{ine22c2}
P(\max\limits_{i=1,\cdots,p}\v \hat{\bw}_{i1} - \bw_{i1} \v \geq t)
&= P(\max\limits_{i=1,\cdots,p} \v \bSXni^{-1}\mathbf{f}_i -\bsig_i^{-1} \mathbf{f}_i\v \geq t)
 \leq \sum_{i=1}^p P(\v \bSXni^{-1}\mathbf{f}_i - \bsig_i^{-1} \mathbf{f}_i \v \geq t) \notag \\
 &\leq \sum_{i=1}^p P(\v \bSXni^{-1} - \bsig_i^{-1}\v \geq t) \leq \sum\limits_{i=1}^p P(\v \bSXni - \bsig_i\v \geq t), 
 \end{align}
notice that the last inequality is due to (\ref{DisInv}) and for simplicity in notation we still use $t$, not ${t}/{C}$. Using the decomposition $\bSXni = n\bVXni/(n-1) + n\bar{\bX}_{\bB_i}\bar{\bX}_{\bB_i}^T/(n-1)$, we have
 \begin{align}
P(\max\limits_{i=1,\cdots,p}\v \hat{\bw}_{i1} - \bw_{i1} \v \geq t) &\leq \sum\limits_{i=1}^p P(\v \bVXni - \bsig_i\v \geq \frac{t}{2}) + \sum\limits_{i=1}^p P(\v \bar{\bX}_{\bB_i}\bar{\bX}_{\bB_i}^T \v \geq \frac{t}{2})\notag \\
 &\leq \sum\limits_{i=1}^p P(s_0\v \bVXni - \bsig_i\v_\infty \geq \frac{t}{2}) 
+ \sum\limits_{i=1}^p P(s_0 \max\limits_{1 \leq j \leq s_0} \l \bar{X}_{\bB_i,j} \l^2 \geq \frac{t}{2}).\label{theo2ab}
\end{align}

For the first term on the right hand side of the inequality, applying Lemma A.3 in Bickel and Levina (2008), we have
\begin{align*}
\sum\limits_{i=1}^p P(\v \bVXni - \bsig_i \v \geq \frac{t}{2})   
\leq \sum\limits_{i=1}^p P(s_0\v \bVXni - \bsig_i \v_\infty \geq \frac{t}{2}) \leq \sum\limits_{i=1}^p s_0^2\exp\big\{{-nt^2\gamma(C_0,\lambda)}/{(4s_0^2)}\big\},
\end{align*}  
for $\l t \l \leq \lambda \equiv \lambda(C_0)$ and some $\gamma(C_0,\lambda) > 0$ that depends only on $C_0$ and $\lambda$.

Choose $t = 2M\sqrt{{s_0^2 \log (s_0^2p)}/{n}}$, for $M > 0$ sufficient large, then we have
\begin{align}\label{theo2a}
\sum\limits_{i=1}^p P(\v  \bVXni - \bsig_i \v \geq \frac{t}{2})
  \leq \sum\limits_{i=1}^p s_0^2\exp\{{-4M^2s_0^2\log(s_0^2p) \gamma(C_0,\lambda)}/{(4s_0^2)}\} 
   \leq \sum\limits_{i=1}^p s_0^2\exp\{ \log (s_0^2p)^{-M^2 \gamma(C_0,\lambda)}\} = o(1).
\end{align}

For the second term, the following inequality holds. For some constant $C > 0$,
 \begin{align} \label{theo2b}
 \sum\limits_{i=1}^p P(\v \bar{\bX}_{\bB_i} \bar{\bX}_{\bB_i}^T \v \geq \frac{t}{2}) 
 &\leq  \sum\limits_{i=1}^p P(s_0 \max\limits_{1 \leq j \leq s_0} \l \bar{X}_{\bB_i,j} \l^2 \geq \frac{t}{2})\leq  \sum\limits_{i=1}^p P\{ \max\limits_{1 \leq j \leq s_0} 
 \l \bar{X}_{\bB_i,j}  \l \geq \sqrt{{t}/(2s_0)}\} \notag \\
 &\leq  \sum\limits_{i=1}^p \sum\limits_{j=1}^{s_0} 
 P\{  \l \bar{X}_{\bB_i,j} \l \geq \sqrt{{t}/(2s_0)}\} \leq 2\sum\limits_{i=1}^p \sum\limits_{j=1}^{s_0} \exp\{- C n^2t/(2s_0)\}.
 \end{align}

Choose $t = M_1{s_0\log (p s_0)}/{n^2}$ for $M_1$ large enough, then $\sum\limits_{i=1}^p P(\v  \bar{\bX}_{\bB_i}\bar{\bX}_{\bB_i}^T \v \geq {t}/{2}) \leq (ps_0)^{1 -C} = o(1),$ for some large $C > 0$.
It can be seen that ${s_0\log (p s_0)}/{n^2}  = o(s_0 \sqrt{\log p/n})$ and $ 2M\sqrt{{s_0^2 \log (ps_0^2)}/{n}} = M_1 s_0 \sqrt{ \log p/n}$, for some $M_1 > 0 $. These facts together with (\ref{theo2ab}), (\ref{theo2a}), and (\ref{theo2b}) we conclude the result $\max\limits_{i=1,\cdots,p}\v \hat{\bw}_{i1} - \bw_{i1} \v = O_P(s_0\sqrt{{\log p}/{n}})$.  \hfill  $\square$

\smallskip

\noindent \textit{Proof of Theorem 3:}
If $\|\bOmg\|_1<\infty$, using Theorem 2.1. in Ju\'{a}rez-Ruiz et al. (2016), we can show that there exists some constants $C > 0 $, such that $ \v \bOmg_i \v_1 < C$, for any $i = 1,2,\cdots, p$.
We have
\begin{align*}
P(\max\limits_{i,j =1,\cdots,p} \l \hat{\omega}_{ij} - \omega_{ij} \l \geq t) \leq \sum\limits_{j =1}^p P( \l \hat{\bw}_{j} - \bw_{j} \l_\infty \geq t)\leq \sum\limits_{j =1}^p P( \l \hat{\bw}_{j1} - \bw_{j1} \l_\infty \geq t) 
= \sum\limits_{j =1}^p P( \l \bSXnj^{-1}\mathbf{f}_j -  \bsig_j^{-1} \mathbf{f}_j \l_\infty \geq t).
\end{align*}
Applying Lemma \ref{represent}, we have $\bSXnj^{-1}\mathbf{f}_j -  \bsig_j^{-1} \mathbf{f}_j = - \bOmg_j \bW_j \bOmg_j \mathbf{f}_j - \bOmg_j \bW_j( \bSXnj^{-1} - \bsig_j^{-1} )\mathbf{f}_j.$
Therefore, 
\begin{align*}
P(\max\limits_{i,j =1,\cdots,p} \l \hat{\omega}_{ij} - \omega_{ij} \l \geq t) 
&\leq \sum\limits_{j =1}^p P( \l \bOmg_j \bW_j \bOmg_j \mathbf{f}_j + \bOmg_j \bW_j( \bSXnj^{-1} - \bsig_j^{-1} )\mathbf{f}_j \l_\infty \geq t)\\
&\leq \sum\limits_{j =1}^p P( \l \bOmg_j \bW_j \bOmg_j \mathbf{f}_j  \l_\infty \geq t/2) + \sum\limits_{j =1}^p P( \l  \bOmg_j \bW_j( \bSXnj^{-1} - \bsig_j^{-1} )\mathbf{f}_j \l_\infty \geq t/2).
\end{align*}
Using the inequality $\l AB \l_\infty \leq \mbox{min}\{\l A \l_\infty \v B\v_1, \v A \v_1 \l B\l_\infty\}$, for matrix $A, B$ with the appropriate size, we have:
\begin{align*}
 \sum\limits_{j =1}^p P( \l \bOmg_j \bW_j \bOmg_j \mathbf{f}_j  \l_\infty \geq t/2) &\leq  \sum\limits_{j =1}^p P( \v \bOmg_j \v_1  \l \bW_j \bOmg_j \mathbf{f}_j  \l_\infty \geq t/2)
 \leq \sum\limits_{j =1}^p P( \v \bOmg_j \v_1  \l  \l \bW_j \l_\infty \v \bOmg_j \mathbf{f}_j  \v_1 \geq t/2)\\
 &\leq \sum\limits_{j =1}^p P( \l \bW_j \l_\infty  \geq t/C),\ \mbox{for some}\ C>0.
\end{align*}

In addition,
$ \sum\limits_{j =1}^p P( \l  \bOmg_j \bW_j( \bSXnj^{-1} - \bsig_j^{-1} )\mathbf{f}_j \l_\infty \geq t/2) \leq  \sum\limits_{j =1}^p P( \v  \bOmg_j\v_1 \l \bW_j \l_ \infty \v ( \bSXnj^{-1} - \bsig_j^{-1} )\mathbf{f}_j \v_1 \geq t/2)$
For each $j =1,2,\cdots, p$, applying Theorem \ref{columnrate}, we have
$\v ( \bSXnj^{-1} - \bsig_j^{-1} )\mathbf{f}_j \v_1 \leq \sqrt{s_0 }\v ( \bSXnj^{-1} - \bsig_j^{-1} )\mathbf{f}_j \v = O_P(\sqrt{s_0^2/n})= o_P(1). $
Therefore, $\sum\limits_{j =1}^p P( \l  \bOmg_j \bW_j( \bSXnj^{-1} - \bsig_j^{-1} )\mathbf{f}_j \l_\infty \geq t/2)\leq  \sum\limits_{j =1}^p P(  \l \bW_j \l_ \infty \geq t/C)$, for some $C > 0$.
 Following the same argument as that in the proof of Theorem \ref{uniformrate}, we can show that by choosing $t = M_1 \sqrt{\mbox{log}(p)/n}$, for some large $M_1 > 0$, both of the quantities $\sum\limits_{j =1}^p P( \l \bOmg_j \bW_j \bOmg_j \mathbf{f}_j  \l_\infty \geq t/2)$ and $\sum\limits_{j =1}^p P( \l  \bOmg_j \bW_j( \bSXnj^{-1} - \bsig_j^{-1} )\mathbf{f}_j \l_\infty \geq t/2)$ converge to 0 as n tends to $\infty$.
 In other words, we have shown that $\max\limits_{i,j =1,\cdots,p} \l \hat{\omega}_{ij} - \omega_{ij} \l = O_P(\sqrt{{\log p}/{n}})$. Thereby proving the theorem. \hfill $\square$

The following lemmas \ref{represent}, \ref{resrate}, \ref{variance} are needed to prove the asymptotic normality in the main results section.

 \begin{lemma}\label{represent}  For any matrix $\bA_{m \times s_0}$, where $m$ is a fixed integer and $\bg_i \in \mathbb{R}^{s_0 \times 1}$ is the $i$-th column of $\bI_{s_0}$, we have
$\bA(\bSYn^{-1}-\bOmg_{\bY})\bg_i = -\bA \bOmg_{\bY} \bWYn \bOmg_{\bY} \bg_i - \bA \bOmg_{\bY} \bWYn (\bSYn^{-1}- \bOmg_{\bY}) \bg_i,$
for $i = 1, \cdots, s_0.$
\end{lemma}

\noindent\textit{Proof:} We have
$ \bSYn^{-1} - \bOmg_{\bY} = -\bOmg_{\bY} \bSYn \bSYn^{-1}+ \bOmg_{\bY} \bsig_{\bY} \bSYn^{-1}= -\bOmg_{\bY} \bWYn \bSYn^{-1}
 = -\bOmg_{\bY} \bWYn \bOmg_{\bY} - \bOmg_{\bY} \bWYn (\bSYn^{-1}- \bOmg_{\bY}).$
So $\bA(\bSYn^{-1}-\bOmg_{\bY})\bg_i 
= - \bA \bOmg_{\bY} \bWYn \bOmg_{\bY} \bg_i
 - \bA \bOmg_{\bY} \bWYn (\bSYn^{-1}- \bOmg_{\bY}) \bg_i $\hfill $\square$

\begin{lemma}\label{resrate}  Let $\bm_{s_0 \times 1}= (m_1, \cdots, m_{s_0})^T$. Then
 $\l \bm^T \bOmg_{\bY} \bWYn (\bSYn^{-1}- \bOmg_{\bY}) \bg_i \l=O_P({s_0 \v \bm \v}/{n}).$
\end{lemma}

 
\noindent\textit{Proof:} We have
$ \l \bm^T \bOmg_{\bY} \bWYn (\bSYn^{-1}- \bOmg_{\bY}) \bg_i \l =  \v \bm^T \bOmg_{\bY} \bWYn (\bSYn^{-1}- \bOmg_{\bY}) \bg_i \v 
 \leq \v \bm  \v \ \v \bOmg_{\bY}  \v  \ \v \bWYn \v  \  \v \bSYn^{-1}- \bOmg_{\bY} \v \ \v \bg_i\v. $
Notice that $\v \bOmg_{\bY}  \v $ = $O(1)$, $ \v \bWYn \v $ =  $O_P(\sqrt{{{s_0}/{{n}}}})$ due to (\ref{consistent3} ).
In addition, Theorem \ref{columnrate} gives us $\v \bSYn^{-1}- \bOmg_{\bY} \v  = O_P(\sqrt{{s_0}/{n}}).$ 
So $\l\bm^T \bOmg_{\bY} \bWYn (\bSYn^{-1}- \bOmg_{\bY}) \bg_i \l=O_P({s_0 \v \bm \v}/{n})$ \hfill $\square$

\begin{lemma}\label{variance}  For any $\ba_{s_0 \times 1} = (a_1,\cdots,a_{s_0})^T$,
  $\bb_{s_0 \times 1} = (b_1,\cdots,b_{s_0})$, where $a_i, b_j $ are bounded for every $1 \leq i, j \leq s_0$, we have $ \Var(\ba^T \bY_1\bY_1^T \bb) =  \ba^T \bsig_{\bY} \ba \bb^T \bsig_{\bY} \bb +  {(\ba^T  \bsig_{\bY} \bb)}^2. $
  In particular, if $\ba^T = \bm^T \Omega $, $\bb = \bOmg_{\bY} \bg_i$, where 
  $\bm_{s_0 \times 1} = (m_1,\cdots,m_{s_0})^T$, 
  $m_i = O(1), \forall i = 1,\cdots,s_0$, then 
 $ \Var(\ba^T \bY_1 \bY_1^T \bb) 
 = \Var(\bm^T \bOmg_{\bY} \bY_1 \bY_1^T \bOmg_{\bY} \bg_i) 
 \asymp {O(\bm^T \bOmg_{\bY} \bm)}.$
 \end{lemma}

 
\noindent\textit{Proof:}
We have $\Var(\ba^T \bY_1 \bY_1^T \bb) = \E(\ba^T \bY_1 \bY_1^T \bb)^2  - \{\E(\ba^T \bY_1 \bY_1^T \bb)\}^2 = \E(\ba^T \bY_1 \bY_1^T \bb)^2 - (\ba^T \bsig_{\bY} \bb)^2.$
Let $\bsig_{\bY} = \bgam^T \bgam$, 
we have $\bgam^T \bU \overset{d}= \bY_1 \sim N(0,\bsig_{\bY})$ for $\bU \sim N(0, \bI_{s_0}).$
Now we calculate $\E(\ba^T \bY_1 \bY_1^T \bb)^2$, we have
$\E(\ba^T \bY_1 \bY_1^T \bb)^2 = \E(\ba^T \bY_1 \bY_1^T \bb \bb^T \bY_1 \bY_1^T \ba ) = \E( \bY_1^T \bb \bb^T \bY_1 \bY_1^T \ba \ba^T \bY_1 )
 = \E( \bU^T \bgam \bb \bb^T \bgam^T \bU \bU^T \bgam \ba \ba^T  \bgam^T \bU) =  \E( \bU^T \bN_1 \bU \bU^T \bN_2 \bU)$
where $\bN_1 = \bgam \bb \bb^T \bgam^T$ and $\bN_2 = \bgam \ba \ba^T  \bgam^T$.

Applying the moment formula for quadratic form on page 188 of Ullah (2004), we get $ \E( \bU^T \bN_1 \bU \bU^T \bN_2 \bU) = \tr (\bN_1)\tr (\bN_2) + 2\tr(\bN_1 \bN_2).$
Moreover, we have $\tr(\bN_1) = \tr (\bgam \bb \bb^T \bgam^T) = \tr (\bb^T \bgam^T \bgam \bb) = \tr (\bb^T \bsig_{\bY} \bb) = \bb^T \bsig_{\bY} \bb$,
$\tr(\bN_2) = \tr(\bgam \ba \ba^T  \bgam^T)= \tr(\ba^T  \bgam^T \bgam \ba ) = \tr(\ba^T \bsig_{\bY} \ba ) = \ba^T \bsig_{\bY} \ba$ and 
 $\tr(\bN_1 \bN_2) = \tr( \bgam \bb \bb^T \bgam^T \bgam \ba \ba^T  \bgam^T) =  \tr( \bb^T \bgam^T \bgam \ba \ba^T  \bgam^T \bgam \bb ) 
= \tr( \bb^T \bsig_{\bY} \ba \ba^T  \bsig_{\bY} \bb ) = {( \ba^T  \bsig_{\bY} \bb)}^2.$
So  $\E(\ba^T \bY_1 \bY_1^T \bb)^2 $$=  \ba^T \bsig_{\bY} \ba \bb^T \bsig_{\bY} \bb + 2{( \ba^T  \bsig_{\bY} \bb)}^2.$
This gives us $\Var(\ba^T \bY_1 \bY_1^T \bb) = \ba^T \bsig_{\bY} \ba \bb^T \bsig_{\bY} \bb+(\ba^T \bsig_{\bY} \bb)^2.$
For $\ba^T = \bm^T \bOmg_{\bY},$ $ \bb = \bOmg_{\bY} \bg_i$, we have
$\mbox{tr}(\bN_1)= \bb^T \bsig_{\bY} \bb= \bg_i^T \bOmg_{\bY} \bsig_{\bY} \bOmg_{\bY} \bg_i = \omega_{ii},$
$\mbox{tr}(\bN_2)= \ba^T \bsig_{\bY} \ba = (\bm^T \bOmg_{\bY} \bsig_{\bY} \bOmg_{\bY} \bm) = \bm^T \bOmg_{\bY} \bm$ and 
$\mbox{tr}(\bN_1 \bN_2)= {(\ba^T \bsig_{\bY} \bb)}^2 = {(\bm^T \bOmg_{\bY} \bsig_{\bY} \bOmg_{\bY} \bg_i)}^2 = (\bm^T  \bOmg_{\bY} \bg_i)^2.$
So $\Var(\bm^T \bOmg_{\bY} \bY_1 \bY_1^T \bOmg_{\bY} \bg_i) = \omega_{ii} \bm^T \bOmg_{\bY} \bm + (\bm^T \bOmg_{\bY} \bg_i)^2. $
Since $\bm^T \bOmg_{\bY} \bg_i$ = $\sum_{j=1}^{s_0} m_j\omega_{ji}= O(1)$, 
we have $\Var(\bm^T \bOmg_{\bY} \bY_1 \bY_1^T \bOmg_{\bY} \bg_i)\asymp $ $ O(\bm^T \bOmg_{\bY} \bm)$. This proves the lemma \hfill$\square$

The following Lemma \ref{sample.cov.asymp} gives us the asymptotic normality result for linear combinations of the sample covariance matrix. 

 \begin{lemma}
 \label{sample.cov.asymp}
Define $h = \Var( \ba^T \bY_1 \bY_1^T \bb)$. Consider any $\ba_{s_0 \times 1} = (a_1,\cdots,a_{s_0})^T$, $\bb_{s_0 \times 1} = (b_1,\cdots,b_{s_0})'$ 
 where $a_i, b_j $ are bounded for every $1 \leq i, j \leq s_0$. If $\v \bsig_{\bY} \v_1 = O(1)$, then $\sqrt{n}(\ba^T \bSYn \bb - \ba^T \bsig_{\bY} \bb)/\sqrt{h} \stackrel{d}{\to} N(0,1).$
\end{lemma}

\noindent \textit {Proof:}
We note that
\begin{align}
\frac{\sqrt{n}}{\sqrt{h}}(\ba^T \bSYn \bb - \ba^T \bsig_{\bY} \bb) 
&= \frac{\sqrt{n}}{\sqrt{h}}(\frac{1}{n-1}\sum\limits_{l = 1}^n \ba^T \bY_l \bY_l^T  \bb  
- \ba^T \bsig_{\bY} \bb)- \frac{\sqrt{n}}{\sqrt{h}}\frac{n}{n-1}\ba^T \bar{\bY} \ {\bar{\bY}}^T  \bb= I_1 - I_2 .\notag
\end{align} 
First, we show that $I_2$ is ignorable under the condition $s_0 = o(\sqrt{n}).$
Notice that, for each fixed $i$ we have $\sum_{i \neq j}\l \sigma_{ij} \l < \infty $, 
so $\sum_{j =1} ^{s_0} a_{i}b_j{\sigma_{ij}}/{n} \leq H \sum_{i\neq j}{\l \sigma_{ij} \l}/{n} = O({1}/{n}) $
where $H= \max_{i,j=1,\cdots,s_0} \l a_i b_j \l= \max_{i,j=1,\cdots,s_0} \l h_{ij} \l  = O(1)$.
This gives us $\E(\ba^T \bar{\bY}\ {\bar{\bY}}^T  \bb)
= \E(\sum\limits_{i,j=1}^{s_0} a_ib_j \bar{Y}_i  \bar{Y}_j)
 = \sum\limits_{i, j=1}^{s_0} a_i b_j{\sigma_{ij}}/{n} 
 = \sum\limits_{i=1}^{s_0} \sum\limits_{j= 1}^{s_0}
  a_i b_j{\sigma_{ij}}/{n}  = O({s_0}/{n}).$

Since $\sqrt{n}{(\bar{Y}_m, \bar{Y}_n, \bar{Y}_p, \bar{Y}_q)}^T$ are jointly  $N(0, \bsig_0)$, 
where $m, n, p, q$ are different in $1,\cdots,s_0$ and $\bsig_0 = (\sigma_{ij})_{i,j  \in \{m,n,p,q\}}.$
This gives us $\E(\bar{Y}_i \bar{Y}_j \bar{Y}_k \bar{Y}_l)=(\sigma_{ij}\sigma_{kl} + \sigma_{ik}\sigma_{jl}+\sigma_{il}\sigma_{jk})/{n^2}$, 
for any $i, j, k, l \in {1,\cdots, s_0}$. So 
\begin{align*}
n^2 \E(\sum\limits_{i,j=1}^{s_0} a_i b_j\bar{Y}_i \bar{Y}_j)^2 
&= n^2 \E(\sum\limits_{i,j,k,l=1}^{s_0}a_ib_ja_kb_l\bar{Y_i} \bar{Y_j} \bar{Y}_k \bar{Y}_l) 
= \sum\limits_{i,j,k,l=1}^{s_0} a_ib_ja_kb_l(\sigma_{ij}\sigma_{kl} + \sigma_{ik}\sigma_{jl} + \sigma_{il}\sigma_{jk} )\\
&= \sum\limits_{i,j,k,l=1}^{s_0} a_ib_ja_kb_l\sigma_{ij}\sigma_{kl}  
  +  \sum\limits_{i,j,k,l=1}^{s_0}a_ib_ja_kb_l\sigma_{ik}\sigma_{jl}
 +\sum\limits_{i,j,k,l=1}^{s_0} a_ib_ja_kb_l \sigma_{il}\sigma_{jk} = O(s_0^2).
\end{align*}
So, we have $\Var (\sum_{i,j=1}^{s_0}a_ib_j\bar{Y}_i\bar{Y}_j)=O ({s_0^2}/{n^2}).$ By Chebyshev Inequality
 \begin{align*}
 P({\sqrt{n} \ \l \sum_{i,j=1}^{s_0}a_{i} b_{j}\bar{Y}_i\bar{Y}_j \l}
 \geq \epsilon\sqrt{h}) &\leq {n \Var (\sum_{i,j=1}^{s_0}h_{ij}\bar{Y_i} \bar{Y_j})}/{(h \epsilon^2)} = O \{{s_0^2}/(n\epsilon^2)\}\rightarrow 0,
 \end{align*}
  as $n \rightarrow \infty$, $s_0 =o(\sqrt{n})$. So $I_2$ is ignorable as $n \rightarrow \infty$. 

We now turn to $I_1$. We have
 $$\frac{1}{n-1}\sum\limits_{l = 1}^n \ba^T \bY_l \bY_l^T  \bb 
 = \frac{n}{n-1}\frac{1}{n}\sum\limits_{l = 1}^n \ba^T \bY_l \bY_l^T  \bb 
 := \frac{n}{n-1}\frac{1}{n} \sum\limits_{l = 1}^n Z_l,$$ 
 where $Z_l = \ba^T \bY_l \bY_l^T \bb$ and $ Z_1, Z_2,\cdots,Z_n$ are IID. Since $\E(Z_l) = \ba^T \bsig_{\bY} \bb, \Var(Z_l) = h < \infty$, by CLT we get
$$\frac{\sqrt{n}}{\sqrt{h}}(\ba^T \bSYn \bb-\ba^T \bsig_{\bY} \bb) \sim AN(0,1).$$
This finishes the proof of this lemma. \hfill $\square$

\end{document}